# Pareto optimality conditions and duality for vector quadratic fractional optimization problems

W. A. Oliveira[a]*, A. Beato Moreno[b], A. C. Moretti[a] and L. L. Salles Neto[c]

[a]*School of Applied Science, University of Campinas, Limeira, SP, Brazil*

[b]*Department of Statistics and Operations Research, College of Mathematics, University of Sevilla, 41012, Sevilla, Spain*

[c]*Department of Science and Technology, Federal University of São Paulo, São José dos Campos, SP, Brazil*

*E-mail: washington.oliveira@fca.unicamp.br [Oliveira]; beato@us.es [Beato Moreno]; moretti@ime.unicamp.br [Moretti]; luiz.leduino@unifesp.br [Salles Neto]*

**Abstract**

One of the most important optimality conditions to aid to solve a vector optimization problem is the first-order necessary optimality condition that generalizes the Karush-Kuhn-Tucker condition. However, to obtain the sufficient optimality conditions, it is necessary to impose additional assumptions on the objective functions and in the constraint set. The present work is concerned with the constrained vector quadratic fractional optimization problem. It shows that sufficient Pareto optimality conditions and the main duality theorems can be established without the assumption of generalized convexity in the objective functions, by considering some assumptions on a linear combination of Hessian matrices instead. The main aspect of this contribution is the development of Pareto optimality conditions based on a similar second-order sufficient condition for problems with convex constraints, without convexity assumptions on the objective functions. These conditions might be useful to determine termination criteria in the development of algorithms.

## 1 Introduction

There are many contributions, concepts, and definitions that characterize and give the Pareto optimality conditions for solutions of a vector optimization problem (see, for instance [9, 28]). One of the most important is the first-order necessary optimality condition that generalizes the Karush-Kuhn-Tucker (KKT) condition. However, to obtain the sufficient optimality conditions, it is necessary to impose additional assumptions (like convexity and its generalizations) in the objective functions and in the constraint set.

*Corresponding author

In this paper, we deal with a particular case of vector optimization problem (VOP), where each objective function consists of a ratio of two quadratic functions. Without generalized convexity assumptions in the objective functions, but by imposing some additional assumptions on a linear combination of Hessian matrices, Pareto optimality conditions are obtained and duality theorems are established. Let us consider the following vector quadratic fractional optimization problem:

$$\begin{array}{lll} \textbf{(VQFP)} & \textbf{Minimize} & \frac{f(x)}{g(x)} = \left(\frac{f_1(x)}{g_1(x)}, \ldots, \frac{f_m(x)}{g_m(x)}\right) \\ & \text{subject to} & h_j(x) \leqq 0 \quad j \in J, \\ & & x \in \Omega, \end{array}$$

where $\Omega \subseteq \mathbb{R}^n$ is an open set, $f_i$, $g_i$, $i \in I \equiv \{1, \ldots, m\}$, and $h_j$, $j \in J \equiv \{1, \ldots, \ell\}$, are continuously differenciable real-valued functions defined on $\Omega$. In addition, we assume that $f_i$, $g_i$, $i \in I$, are quadratic functions and $g_i(x) > 0$ for $x \in \Omega$ and $i \in I$. We denote by $S$ the *feasible set* of elements $x \in \Omega$ satisfying $h_j(x) \leqq 0$. We say that $x$ is a *feasible point* if $x \in S$. The value $\frac{f_i(x)}{g_i(x)}$ is the result of the $i^{th}$ objective function if the decision maker chooses the action $x \in S$.

Fractional optimization problems arise frequently in decision making applications, including science management, portfolio selection, cutting and stock, game theory, in the optimization of the ratio performance/cost, or profit/ investment, or cost/time and so on.

There are many contributions dealing with the scalar (single-objective) fractional optimization problem (FP) and vector fractional optimization problem (VFP). In most of them, using convexity or their generalizations, optimality conditions in the KKT sense, and the main duality theorems for optimal points are obtained. With a parametric approach, which transforms the original problem in a simpler associated problem, Dinkelbach [12], Jagannathan [16] and Antczak [1] established optimality conditions, presented algorithms and applied their approaches in a example (FP) consisting of quadratic functions. Using some known generalized convexity, Antczak [1], Khan and Hanson [19], Reddy and Mukherjee [32], Jeyakumar [17], Liang et al. [24] established optimality conditions and theorems that relate the pair primal-dual of problem (FP). In Craven [10] and Weir [37], other results for the scalar optimization (FP) can be found.

Further, Liang et al. [25] extended their approach to the vector optimization case (VFP) considering the type duals of Mond and Weir [29], Schaible [35] and Bector [5]. Considering the parametric approach of Dinkelbach [12], Jagannathan [16], Bector et al. [6] and two classes of generalized convexity, Osuna-Gómez et al. [30] established weak Pareto optimality conditions and the main duality theorems for the differenciable vector optimization case (VFP). Santos et al. [34] deepened these results to the more general non-differenciable case (VFP). Jeyakumar and Mond [18] used generalized convexity to study the problem (VFP).

Few studies are found involving quadratic functions at both the numerator and denominator of the ratio objective function. Most of them involve the mixing of linear and quadratic functions. The most similar approaches to the scalar quadratic fractional optimization problem (QFP) were considered in [8, 11, 14, 26, 36]. On the other hand, Benson [7] considered a pure (QFP) consisting of convex functions, where some theoretical properties and optimality conditions are developed, and an algorithm and its convergence properties are presented.

The closest approaches to the vector optimization case (VQFP) were considered in [2, 3, 4, 20, 21, 22, 23, 33]. Using an iterative computational test, Beato et al. [3, 4] characterized the Pareto optimal point for the problem (VQFP), consisting of a linear and quadratic functions, and some theoretical results were obtained by using the function linearization technique of Bector et al. [6]. Arévalo and Zapata [2], Konno and Inori [20], Rhode and Weber [33] analyzed the portfolio selection problem. Kornbluth and Steuer [23] used an adapted Simplex method in the problem (VFP) consisting of linear functions. Korhonen and Yu [21, 22] proposed an iterative computational method for solving the problem (VQFP), consisting of linear and quadratic functions, based on search directions and weighted sums.

The approach taken in this work is different from the previous ones. The main aspect of this contribution is the development of Pareto optimality conditions for a particular vector optimization problem based on a similar second-order sufficient condition for Pareto optimality for problems with convex constraints without the hypothesis of convexity in the objective functions. These conditions might be useful to determine termination criteria in the development of algorithms, and new extensions can be established from these, where more general vector optimization problems in which algorithms based on quadratic approximations are used locally.

This paper is organized as follows. We start by defining some notations and basic properties in Section 2. In Section 3, the sufficient Pareto optimality conditions are established. In Section 4, the relationship among the associated problems is presented and duality theorems are established. Finally, comments and concluding remarks are presented in Section 5.

# 2 Preliminaries

Let $\mathbb{R}_+$ denote the nonnegative real numbers and $x^T$ denote the transpose of the vector $x \in \mathbb{R}^n$. Furthermore, we will adopt the following conventions for inequalities among vectors. If $x = (x_1, \ldots, x_m)^T \in \mathbb{R}^m$ and $y = (y_1, \ldots, y_m)^T \in \mathbb{R}^m$, then

$$
\begin{aligned}
x = y \quad & \text{if and only if} \quad x_i = y_i, \quad \forall i \in I; \\
x < y \quad & \text{if and only if} \quad x_i < y_i, \quad \forall i \in I; \\
x \leqq y \quad & \text{if and only if} \quad x_i \leqq y_i, \quad \forall i \in I; \\
x \leq y \quad & \text{if and only if} \quad x \leqq y \text{ and } x \neq y.
\end{aligned}
$$

Similarly we consider the equivalent convention for inequalities $>, \geqq$ and $\geq$.

Different optimality definitions for the problem (VQFP) are referred as Pareto optimal solutions [31], two of which are defined as follows.

**Definition 1** A feasible point $x^*$ is said to be a *Pareto optimal solution* of (VQFP), if there does not exist another $x \in S$ such that $\frac{f(x)}{g(x)} \leq \frac{f(x^*)}{g(x^*)}$.

**Definition 2** A feasible point $x^*$ is said to be a *weakly Pareto optimal solution* of (VQFP), if there does not exist another $x \in S$ such that $\frac{f(x)}{g(x)} < \frac{f(x^*)}{g(x^*)}$.

Hypotheses of convexity or generalized convexity on the objective functions will be avoided in this work, but we will use such hypotheses in the constraint set. We recall the definition of convexity, where $\nabla f(x)$ denotes the gradient of the function $f : \mathbb{R}^n \to \mathbb{R}$ at the point $x$.

**Definition 3** Let $f : \Omega \subseteq \mathbb{R}^n \to \mathbb{R}$ be a function defined on an open convex set $\Omega$ and differenciable at $x^* \in \Omega$. $f$ is called *convex* at $x^*$ if for all $x \in \Omega$, $f(x) - f(x^*) \geqq \nabla f(x^*)^T (x - x^*)$. When $f$ is convex on the set $\Omega$, we simply say that $f$ is convex.

Maeda [27] used the generalized Guignard constraint qualification (GGCQ) [15] to derive the following necessary Pareto optimality conditions for the problem (VOP) in the KKT sense. Assuming differentiability of the objective and the constraint functions, Maeda guarantees the existence of Lagrange multipliers, all strictly positive, associated with the objective functions.

**Lemma 2.1** *(Maeda [27]) Let $x^*$ be a Pareto optimal solution of (VQFP). Suppose that (GGCQ) holds at $x^*$, then there exist vectors $\tau \in \mathbb{R}^m$, $\lambda \in \mathbb{R}^\ell$ such that*

$$\begin{aligned}
&\sum_{i=1}^{m} \tau_i \nabla \frac{f_i(x^*)}{g_i(x^*)} + \sum_{j=1}^{\ell} \lambda_j \nabla h_j(x^*) = 0,\\
&\sum_{j=1}^{\ell} \lambda_j h_j(x^*) = 0,\\
&\tau > 0, \quad \lambda \geqq 0.
\end{aligned}$$

For each $i \in I$ and $x \in \mathbb{R}^n$, we consider the objective functions defined as $f_i(x) = x^T A_i x + a_i^T x + \bar{a}_i$ and $g_i(x) = x^T B_i x + b_i^T x + \bar{b}_i$, where $A_i$, $B_i \in \mathbb{R}^{n\times n}$, $A_i$ symmetric, $B_i$ symmetric and positive semidefinite, $a_i$, $b_i \in \mathbb{R}^n$ and $\bar{a}_i$, $\bar{b}_i \in \mathbb{R}$, with $\bar{b}_i > -({w^i}^T B_i w^i + b_i^T w^i)$, where $w^i$ is the solution of the system $2B_i x + b_i = 0$, that is, $w^i$ is the point in which the function $x^T B_i x + b_i^T x$ reaches its minimum and this ensures that $g_i(x) > 0$, $\forall x \in \mathbb{R}^n$. We cannot consider the cases where $2B_i x + b_i = 0$ has no solution.

# 3 Sufficient optimality conditions

Without assumptions of generalized convexity, but imposing some additional assumptions on a linear combination of Hessian matrices of the objective functions $f_i$ and $g_i$, $i \in I$, we provide in the next theorem a sufficient condition that guarantees that a feasible point of (VQFP) is Pareto optimal point. Similar to a second-order sufficient condition for Pareto optimality, this condition explores the intrinsic characteristics of the problem (VQFP).

We assume, unlike the objective functions, that each $h_j$ is convex. Also, given $x^* \in S$, for each $i \in I$ we define the scalar functions $u_i : S \times S \to \mathbb{R}_+ \setminus \{0\}$ and $s_i : S \times S \to \mathbb{R}$ by

$$\begin{aligned}
u_i(x, x^*) &\equiv \frac{g_i(x^*)}{g_i(x)},\\
s_i(x, x^*) &\equiv \frac{1}{g_i(x)} \left\{ (x - x^*)^T \left[ A_i - \frac{f_i(x^*)}{g_i(x^*)} B_i \right] (x - x^*) \right\}.
\end{aligned}$$

**Theorem 3.1** *Let $x^*$ be a feasible point of (VQFP). Suppose that the constraint function $h_j$ is convex for each $j \in J$ and there exist vectors $\tau \in \mathbb{R}^m$, $\lambda \in \mathbb{R}^\ell$, such that*

$$\sum_{i=1}^{m} \tau_i \nabla \frac{f_i(x^*)}{g_i(x^*)} + \sum_{j=1}^{\ell} \lambda_j \nabla h_j(x^*) = 0, \tag{1}$$

$$\sum_{j=1}^{\ell} \lambda_j h_j(x^*) = 0, \tag{2}$$

$$\tau > 0, \quad \lambda \geqq 0. \tag{3}$$

*If for any $x \in S$, we obtain*

$$\sum_{i=1}^{m} \tau_i \frac{s_i(x,x^*)}{u_i(x,x^*)} \geqq 0, \tag{4}$$

*then $x^*$ is a Pareto optimal solution for (VQFP).*

*Proof* Given $x \in S$, we obtain for each $i \in I$

$$\begin{aligned} f_i(x) - f_i(x^*) &= (x-x^*)^T A_i (x-x^*) + \nabla f_i(x^*)^T (x-x^*), \\ g_i(x) - g_i(x^*) &= (x-x^*)^T B_i (x-x^*) + \nabla g_i(x^*)^T (x-x^*), \text{ and} \end{aligned}$$

$$\begin{aligned} & \frac{f_i(x)}{g_i(x)} - \frac{f_i(x^*)}{g_i(x^*)} = \frac{f_i(x)g_i(x^*) - g_i(x)f_i(x^*)}{g_i(x)g_i(x^*)} = \\ = \; & \frac{f_i(x)g_i(x^*) - f_i(x^*)g_i(x^*) + f_i(x^*)g_i(x^*) - g_i(x)f_i(x^*)}{g_i(x)g_i(x^*)} \\ = \; & \frac{g_i(x^*)\left\{f_i(x) - f_i(x^*)\right\} - f_i(x^*)\left\{g_i(x) - g_i(x^*)\right\}}{g_i(x)g_i(x^*)} \\ = \; & \frac{1}{g_i(x)} \left\{ (x-x^*)^T A_i (x-x^*) + \nabla f_i(x^*)^T (x-x^*) \right\} + \\ & \qquad - \frac{f_i(x^*)}{g_i(x)g_i(x^*)} \left\{ (x-x^*)^T B_i (x-x^*) + \nabla g_i(x^*)^T (x-x^*) \right\} \\ = \; & \frac{1}{g_i(x)} \left\{ (x-x^*)^T A_i (x-x^*) \right\} - \frac{f_i(x^*)}{g_i(x)g_i(x^*)} \left\{ (x-x^*)^T B_i (x-x^*) \right\} + \\ & \qquad + \frac{1}{g_i(x)} \left\{ \nabla f_i(x^*)^T (x-x^*) \right\} - \frac{f_i(x^*)}{g_i(x)g_i(x^*)} \left\{ \nabla g_i(x^*)^T (x-x^*) \right\} \\ = \; & \frac{1}{g_i(x)} \left\{ (x-x^*)^T \left[ A_i - \frac{f_i(x^*)}{g_i(x^*)} B_i \right] (x-x^*) \right\} + \\ & \qquad + \frac{g_i(x^*)}{g_i(x)} \left\{ \left[ \frac{\nabla f_i(x^*) g_i(x^*) - \nabla g_i(x^*) f_i(x^*)}{[g_i(x^*)]^2} \right]^T (x-x^*) \right\} \\ = \; & \frac{1}{g_i(x)} \left\{ (x-x^*)^T \left[ A_i - \frac{f_i(x^*)}{g_i(x^*)} B_i \right] (x-x^*) \right\} + \frac{g_i(x^*)}{g_i(x)} \left\{ \left[ \nabla \frac{f_i(x^*)}{g_i(x^*)} \right]^T (x-x^*) \right\}. \end{aligned}$$

Thus, each function $\frac{f_i}{g_i}$ satisfies

$$\frac{f_i(x)}{g_i(x)} - \frac{f_i(x^*)}{g_i(x^*)} = u_i(x,x^*)\left[\nabla \frac{f_i(x^*)}{g_i(x^*)}\right]^T (x-x^*) + s_i(x,x^*), \quad \forall x \in S. \tag{5}$$

Suppose that $x^*$ is not a Pareto optimal solution of (VQFP). Then there exists another point $x \in S$ such that

$$\frac{f(x)}{g(x)} \leq \frac{f(x^*)}{g(x^*)}. \tag{6}$$

Since $u_i(x,x^*) > 0$, $i \in I$, from Equation (5) we obtain

$$\frac{1}{u_i(x,x^*)}\left(\frac{f_i(x)}{g_i(x)} - \frac{f_i(x^*)}{g_i(x^*)}\right) = \left[\nabla \frac{f_i(x^*)}{g_i(x^*)}\right]^T (x-x^*) + \frac{s_i(x,x^*)}{u_i(x,x^*)}, \quad i \in I.$$

From (6), we have

$$\frac{f(x)}{g(x)} - \frac{f(x^*)}{g(x^*)} \leq 0,$$

and we obtain $m$ inequalities

$$\left[\nabla \frac{f_i(x^*)}{g_i(x^*)}\right]^T (x-x^*) + \frac{s_i(x,x^*)}{u_i(x,x^*)} \leqq 0, \quad i \in I,$$

with at least one strict inequality. Multiplying the $m$ inequalities above by their respective $\tau_i > 0$, $i \in I$, and summing all the products, we obtain

$$\sum_{i=1}^{m} \tau_i \left[\nabla \frac{f_i(x^*)}{g_i(x^*)}\right]^T (x-x^*) + \sum_{i=1}^{m} \tau_i \frac{s_i(x,x^*)}{u_i(x,x^*)} < 0.$$

Then, we have

$$\left[\sum_{i=1}^{m} \tau_i \nabla \frac{f_i(x^*)}{g_i(x^*)}\right]^T (x-x^*) + \sum_{i=1}^{m} \tau_i \frac{s_i(x,x^*)}{u_i(x,x^*)} < 0. \tag{7}$$

Substituting (1) into (7), we get

$$\left[-\sum_{j=1}^{\ell} \lambda_j \nabla h_j(x^*)\right]^T (x-x^*) + \sum_{i=1}^{m} \tau_i \frac{s_i(x,x^*)}{u_i(x,x^*)} < 0. \tag{8}$$

Using (4) and (8), we obtain

$$0 \leqq \sum_{i=1}^{m} \tau_i \frac{s_i(x,x^*)}{u_i(x,x^*)} < \left[\sum_{j=1}^{\ell} \lambda_j \nabla h_j(x^*)\right]^T (x-x^*).$$

That is,

$$\left[\sum_{j=1}^{\ell} \lambda_j \nabla h_j(x^*)\right]^T (x-x^*) > 0. \tag{9}$$

On the other hand, by convexity of $h_j$, we have for each $j \in J$,

$$h_j(x) - h_j(x^*) \geqq \nabla h_j(x^*)^T (x-x^*).$$

Since $\lambda_j \geqq 0$, $j \in J$, we have

$$\sum_{j=1}^{\ell} \lambda_j \left(h_j(x) - h_j(x^*)\right) \geqq \sum_{j=1}^{\ell} \lambda_j \nabla h_j(x^*)^T (x-x^*) = \left[\sum_{j=1}^{\ell} \lambda_j \nabla h_j(x^*)\right]^T (x-x^*).$$

However, since $x$ is feasible point, condition (2) and $\lambda_j \geqq 0$, $j \in J$, imply that

$$\sum_{j=1}^{\ell} \lambda_j \left(h_j(x) - h_j(x^*)\right) \leqq 0.$$

We conclude that

$$\left[\sum_{j=1}^{\ell} \lambda_j \nabla h_j(x^*)\right]^T (x-x^*) \leqq 0,$$

which contradicts (9). Therefore $x^*$ is a Pareto optimal solution for (VQFP). ∎

The expression $\frac{f_i(x)}{g_i(x)} - \frac{f_i(x^*)}{g_i(x^*)}$ in Theorem 3.1 is manipulated in a similar manner in [18, 19, 24, 25, 32], however some generalized convexity on the functions $f_i$ and $g_i$ are imposed. In most of them, for each $i \in I$ and $x \in S$, the hypothesis $f_i(x) \geqq 0$, $g_i(x) > 0$ and $f_i$, $-g_i$ satisfy some generalized convexity. This is not the purpose of this work, but the constraint functions can be assumed in a more general class of convex functions, for example, the generalized convexity of Liang et al. [24] can be used.

In the following, the Pareto optimal solution set is denoted by $\mathit{Eff}(VQFP)$.

**Corollary 3.2** *Let $x^*$ be a feasible point of (VQFP). Suppose that the constraint function $h_j$ is convex for each $j \in J$, and there exist vectors $\tau \in \mathbb{R}^m$, $\lambda \in \mathbb{R}^\ell$, such that (1), (2) and (3) are valid. If $\left[A_i - \frac{f_i(x^*)}{g_i(x^*)} B_i\right]$ are positive semidefinite matrices for each $i \in I$, then $x^* \in \mathit{Eff}(VQFP)$.*

*Proof* By hypothesis, given $x \in S$ and $i \in I$, we obtain

$$\begin{aligned}
(x-x^*)^T \left[A_i - \frac{f_i(x^*)}{g_i(x^*)} B_i\right] (x-x^*) &\geqq 0 \Longrightarrow \\
\tau_i \frac{1}{g_i(x^*)} \left\{ (x-x^*)^T \left[A_i - \frac{f_i(x^*)}{g_i(x^*)} B_i\right] (x-x^*) \right\} &\geqq 0 \Longrightarrow \\
\sum_{i=1}^{m} \tau_i \frac{g_i(x)}{g_i(x^*)} \frac{1}{g_i(x)} \left\{ (x-x^*)^T \left[A_i - \frac{f_i(x^*)}{g_i(x^*)} B_i\right] (x-x^*) \right\} &\geqq 0 \Longrightarrow \\
\Longrightarrow \sum_{i=1}^{m} \tau_i \frac{s_i(x,x^*)}{u_i(x,x^*)} &\geqq 0.
\end{aligned}$$

Therefore, the inequality (4) is valid and the result follows from Theorem 3.1. ■

To ensure that inequality (4) is valid, we start exploring the features of the Hessian matrices of the objective functions of (VQFP).

Negative values can occur in each term $\tau_i \frac{s_i(x,x^*)}{u_i(x,x^*)}$ of the sum $\sum_{i=1}^{m} \tau_i \frac{s_i(x,x^*)}{u_i(x,x^*)}$, which depends on each matrix $\left[A_i - \frac{f_i(x^*)}{g_i(x^*)} B_i\right]$, $i \in I$, and the vector $(x - x^*)$. Let us check new conditions for which (4) is satisfied, that is, we want to ensure the result of Theorem 3.1 by analysing the function

$$\begin{aligned} Z(x,x^*) \equiv \sum_{i=1}^{m} \tau_i \frac{s_i(x,x^*)}{u_i(x,x^*)} &= \sum_{i=1}^{m} \tau_i \frac{g_i(x)}{g_i(x^*)} \frac{1}{g_i(x)} \left\{ (x-x^*)^T \left[A_i - \frac{f_i(x^*)}{g_i(x^*)} B_i\right] (x-x^*) \right\} \\ &= (x-x^*)^T \left\{ \sum_{i=1}^{m} \left[ \frac{\tau_i}{g_i(x^*)} A_i - \frac{\tau_i f_i(x^*)}{[g_i(x^*)]^2} B_i \right] \right\} (x-x^*). \end{aligned}$$

Note that $Z(\cdot,x^*)$ is a quadratic function without the linear part, thus in (VQFP) we obtain $Z(x,x^*) \geqq 0$ on $S$ if and only if $\min_{x \in S} Z(x,x^*) \geqq 0$, that is, we can use the classical results on quadratic optimization to check if $\min_{x \in S} Z(x,x^*) \geqq 0$. The next corollary follows immediately from Theorem 3.1.

**Corollary 3.3** *Let $x^*$ be a feasible point of (VQFP). Suppose that the constraint function $h_j$ is convex for each $j \in J$ and there exist vectors $\tau \in \mathbb{R}^m$, $\lambda \in \mathbb{R}^\ell$, such that (1), (2) and (3) are valid. If $\min_{x \in S} Z(x,x^*) \geqq 0$, then $x^* \in Eff(VQFP)$.* ■

Using the previous results to check whether a feasible point is a Pareto optimal solution of (VQFP), we propose the following computational test method.

---

*Pareto optimality test*

*Step 1.* Given $x^* \in S$. Find the vectors $\tau > 0$ and $\lambda \geqq 0$ such that (1) and (2) are valid. If the vectors $\tau$ and $\lambda$ do not exist, then $x^* \notin Eff(VQFP)$.

*Step 2.* Otherwise, solve $Z(\bar{x},x^*) = \min_{x \in S} Z(x,x^*)$. If $Z(\bar{x},x^*) \geqq 0$, we say that $x^*$ has passed the Pareto optimality test and $x^* \in Eff(VQFP)$.

---

Pareto optimality test starts with a feasible point, then it seeks to solve a system of linear equations containing $m + \ell$ unknowns, $\tau$ and $\lambda$, the inequalities $\tau > 0$, $\lambda \geqq 0$, and two equalities (1) and (2). If this system has no solution, then the point $x^*$ does not satisfy the first-order necessary condition for Pareto optimality, so the method terminates concluding that $x^* \notin Eff(VQFP)$. Otherwise, in Step 2, a quadratic optimization problem on $S$ should be solved. If the minimum of the quadratic problem is non-negative, then the procedure ends, concluding that $x^* \in Eff(VQFP)$. Otherwise, we say that $x^*$ has not passed the Pareto optimality test. Its complexity lies in solving a system of linear inequalities plus a quadratic optimization problem.

The next results, which addresses a linear combination of the Hessian matrices, can be used to develop a computational search method.

Looking at the previous Pareto optimality test, if the fixed point $x^*$ is assumed to be a variable $y$, then the linear system in Step 1 becomes a nonlinear system for the variables $\tau > 0$, $\lambda \geqq 0$, $y \in S$. And the quadratic optimization problem in Step 2 becomes a quadratic optimization problem of the type $\min\limits_{x,y \in S} Z(x,y)$. This raises considerable difficulties. In order to reduce these difficulties, we further explore the characteristics of the matrix

$$\widehat{F}(y^*) \equiv \sum_{i=1}^{m} \left[ \frac{\tau_i}{g_i(y^*)} A_i - \frac{\tau_i f_i(y^*)}{[g_i(y^*)]^2} B_i \right]. \tag{10}$$

One possibility is to search for points $y^*$ such that $\widehat{F}(y^*)$ becomes positive semidefinite. In this case $Z(x,y^*) = (x-y^*)^T \widehat{F}(y^*)(x-y^*) \geqq 0$ depends only on $y^* \in S$.

Consider a fixed point $x^*$, the next theorem takes advantage of the symmetry and diagonalizations of the matrices $A_i$ and $B_i$, $i \in I$, to give sufficient Pareto optimality conditions for a feasible point of (VQFP). Consider the usual inner product $\langle \cdot, \cdot \rangle$ in $\mathbb{R}^n$.

**Theorem 3.4** *Let $x^*$ be a feasible point of (VQFP). Suppose that the constraint function $h_j$ is convex for each $j \in J$ and there exist vectors $\tau \in \mathbb{R}^m$, $\lambda \in \mathbb{R}^\ell$, such that (1), (2) and (3) are valid. Consider also, for each $i \in I$ and $k \in K \equiv \{1,\ldots,n\}$, the following functions*

$$\gamma_i^k(x,x^*,\tau) \equiv \frac{\tau_i}{g_i(x^*)} \left\langle x-x^*, p_i^k \right\rangle^2 \quad \text{and} \quad \eta_i^k(x,x^*,\tau) \equiv \frac{\tau_i f_i(x^*)}{[g_i(x^*)]^2} \left\langle x-x^*, q_i^k \right\rangle^2,$$

*where $p_i^k$ and $q_i^k$ are the columns of orthogonal matrices $P_i$ and $Q_i$, constructed from the normalized eigenvectors of the matrices $A_i$ and $B_i$, respectively. If for all $x \in S$ the following inequalities*

$$\mu_k^{A_i} \gamma_i^k(x,x^*,\tau) \geqq \mu_k^{B_i} \eta_i^k(x,x^*,\tau), \ \forall i \in I \text{ and } \ \forall k \in K, \tag{11}$$

*are valid, where $\mu_k^{A_i}$ and $\mu_k^{B_i}$ are the eigenvalues of the matrices $A_i$ and $B_i$ associated with the eigenvectors $p_i^k$ and $q_i^k$, respectively. Then $x^* \in Eff(VQFP)$.*

*Proof* The matrices $A_i$ and $B_i$, $i \in I$, are diagonalizable and can be rewritten as $A_i = P_i D_{A_i} P_i^T = \sum\limits_{k=1}^{n} \mu_k^{A_i} p_i^k p_i^{kT}$ and $B_i = Q_i D_{B_i} Q_i^T = \sum\limits_{k=1}^{n} \mu_k^{B_i} q_i^k q_i^{kT}$, where $D_{A_i}$ and $D_{B_i}$ are diagonal matrices, with their diagonal formed by the eigenvalues $\mu_k^{A_i}$ and $\mu_k^{B_i}$, $k \in K$, of the matrices $A_i$ and $B_i$, respectively. Thus, we obtain

$$\begin{aligned}
&\sum_{i=1}^{m} \tau_i \frac{s_i(x,x^*)}{u_i(x,x^*)} = \sum_{i=1}^{m} \tau_i \frac{g_i(x)}{g_i(x^*)} \frac{1}{g_i(x)} \left\{ (x-x^*)^T \left[ A_i - \frac{f_i(x^*)}{g_i(x^*)} B_i \right] (x-x^*) \right\} = \\
&= \sum_{i=1}^{m} \left\{ (x-x^*)^T \left[ \frac{\tau_i}{g_i(x^*)} A_i - \frac{\tau_i f_i(x^*)}{[g_i(x^*)]^2} B_i \right] (x-x^*) \right\} \\
&= \sum_{i=1}^{m} \left\{ (x-x^*)^T \left[ \frac{\tau_i}{g_i(x^*)} \left( \sum_{k=1}^{n} \mu_k^{A_i} p_i^k p_i^{kT} \right) - \frac{\tau_i f_i(x^*)}{[g_i(x^*)]^2} \left( \sum_{k=1}^{n} \mu_k^{B_i} q_i^k q_i^{kT} \right) \right] (x-x^*) \right\} \\
&= \sum_{i=1}^{m} \left[ \sum_{k=1}^{n} \mu_k^{A_i} \left( \frac{\tau_i}{g_i(x^*)} \left\langle x-x^*, p_i^k \right\rangle^2 \right) - \sum_{k=1}^{n} \mu_k^{B_i} \left( \frac{\tau_i f_i(x^*)}{[g_i(x^*)]^2} \left\langle x-x^*, q_i^k \right\rangle^2 \right) \right]
\end{aligned}$$

$$= \sum_{i=1}^{m} \left[ \sum_{k=1}^{n} \mu_k^{A_i} \gamma_i^k(x,x^*,\tau) - \sum_{k=1}^{n} \mu_k^{B_i} \eta_i^k(x,x^*,\tau) \right]$$
$$= \sum_{i=1}^{m} \sum_{k=1}^{n} \left[ \mu_k^{A_i} \gamma_i^k(x,x^*,\tau) - \mu_k^{B_i} \eta_i^k(x,x^*,\tau) \right].$$

Since for all $x \in S$, we have $\mu_k^{A_i} \gamma_i^k(x,x^*,\tau) \geqq \mu_k^{B_i} \eta_i^k(x,x^*,\tau)$, for all $i \in I$ and $k \in K$, we conclude that $Z(x,x^*) \geqq 0$. Therefore, the inequality (4) is valid and the result follows from Theorem 3.1. ■

Theorem 3.4 is not simple to use since (11) depends on all points of the feasible set, that is, it depends of the functions $\gamma_i^k(x,x^*,\tau)$, $\eta_i^k(x,x^*,\tau)$, $\forall i \in I$, $\forall k \in K$, and $x \in S$. However, even if for some $i \in I$ and $k \in K$, $\mu_k^{A_i} \gamma_i^k(x,x^*,\tau) < \mu_k^{B_i} \eta_i^k(x,x^*,\tau)$ occurs, the inequality (4) can still be satisfied. In order to obtain (11), we present the next corollary, which follows immediately from the previous theorem.

**Corollary 3.5** *Let $x^*$ be a feasible point of (VQFP). Suppose that the constraint function $h_j$ is convex for each $j \in J$, and there exist vectors $\tau \in \mathbb{R}^m$, $\lambda \in \mathbb{R}^\ell$, such that (1), (2) and (3) are valid. Consider also, for each $i \in I$ and $k \in K$,*

$$a_{i,k}^{+} \equiv \sqrt{\mu_k^{A_i}}\ p_i^k + \sqrt{\mu_k^{B_i} \frac{f_i(x^*)}{g_i(x^*)}}\ q_i^k, \qquad a_{i,k}^{-} \equiv \sqrt{\mu_k^{A_i}}\ p_i^k - \sqrt{\mu_k^{B_i} \frac{f_i(x^*)}{g_i(x^*)}}\ q_i^k, \tag{12}$$

$$\alpha_{i,k} \equiv \left\langle x^*, a_{i,k}^{+} \right\rangle a_{i,k}^{-} + \left\langle x^*, a_{i,k}^{-} \right\rangle a_{i,k}^{+}, \quad \beta_{i,k} \equiv \left\langle x^*, a_{i,k}^{+} \right\rangle \left\langle x^*, a_{i,k}^{-} \right\rangle, \tag{13}$$

$$H_{i,k}(x) \equiv x^T \left[ a_{i,k}^{+} {a_{i,k}^{-}}^{T} \right] x - \alpha_{i,k}^T x + \beta_{i,k}, \tag{14}$$

*where $p_i^k$ and $q_i^k$ are the columns of orthogonal matrices $P_i$ and $Q_i$, constructed from the normalized eigenvectors of the matrices $A_i$, $B_i$, and $\mu_k^{A_i}$, $\mu_k^{B_i}$ are the eigenvalues of the matrices $A_i$ and $B_i$ associated with the eigenvectors $p_i^k$ and $q_i^k$, respectively. If for all $x \in S$, we obtain $H_{i,k}(x) \geqq 0$ for each $i \in I$ and $k \in K$, then $x^* \in Eff(VQFP)$.*

*Proof* According to Theorem 3.4, it is enough to show that for every feasible point, and for all $i \in I$ and $k \in K$, $\mu_k^{A_i} \gamma_i^k(x,x^*,\tau) \geqq \mu_k^{B_i} \eta_i^k(x,x^*,\tau)$ is valid. Given $x \in S$ and a pair $\{i,k\} \in I \times K$, we obtain

$$\mu_k^{A_i} \gamma_i^k(x,x^*,\tau) \geqq \mu_k^{B_i} \eta_i^k(x,x^*,\tau) \iff \frac{\tau_i \mu_k^{A_i}}{g_i(x^*)} \left\langle x - x^*, p_i^k \right\rangle^2 \geqq \frac{\tau_i \mu_k^{B_i} f_i(x^*)}{[g_i(x^*)]^2} \left\langle x - x^*, q_i^k \right\rangle^2$$

$$\iff \mu_k^{A_i} \left\langle x - x^*, p_i^k \right\rangle^2 \geqq \mu_k^{B_i} \frac{f_i(x^*)}{g_i(x^*)} \left\langle x - x^*, q_i^k \right\rangle^2 \iff$$

$$\left( \sqrt{\mu_k^{A_i}}\ \left\langle x - x^*, p_i^k \right\rangle + \sqrt{\mu_k^{B_i} \frac{f_i(x^*)}{g_i(x^*)}}\ \left\langle x - x^*, q_i^k \right\rangle \right) \left( \sqrt{\mu_k^{A_i}}\ \left\langle x - x^*, p_i^k \right\rangle - \sqrt{\mu_k^{B_i} \frac{f_i(x^*)}{g_i(x^*)}}\ \left\langle x - x^*, q_i^k \right\rangle \right) \geqq 0 \iff$$

$$\left( \left\langle x - x^*, \sqrt{\mu_k^{A_i}}\ p_i^k + \sqrt{\mu_k^{B_i} \frac{f_i(x^*)}{g_i(x^*)}}\ q_i^k \right\rangle \right) \left( \left\langle x - x^*, \sqrt{\mu_k^{A_i}}\ p_i^k - \sqrt{\mu_k^{B_i} \frac{f_i(x^*)}{g_i(x^*)}}\ q_i^k \right\rangle \right) \geqq 0 \iff$$

$$\left(\left\langle x-x^*,a_{i,k}^+\right\rangle\right)\left(\left\langle x-x^*,a_{i,k}^-\right\rangle\right) \geqq 0 \iff$$

$$\left(\left\langle x,a_{i,k}^+\right\rangle - \left\langle x^*,a_{i,k}^+\right\rangle\right)\left(\left\langle x,a_{i,k}^-\right\rangle - \left\langle x^*,a_{i,k}^-\right\rangle\right) \geqq 0 \iff$$

$$\left\langle x,a_{i,k}^+\right\rangle\left\langle x,a_{i,k}^-\right\rangle - \left\langle x^*,a_{i,k}^-\right\rangle\left\langle x,a_{i,k}^+\right\rangle - \left\langle x^*,a_{i,k}^+\right\rangle\left\langle x,a_{i,k}^-\right\rangle + \left\langle x^*,a_{i,k}^+\right\rangle\left\langle x^*,a_{i,k}^-\right\rangle \geqq 0 \iff$$

$$\left\langle x,a_{i,k}^+\right\rangle\left\langle x,a_{i,k}^-\right\rangle - \left\langle x,\left(\left\langle x^*,a_{i,k}^+\right\rangle a_{i,k}^- + \left\langle x^*,a_{i,k}^-\right\rangle a_{i,k}^+\right)\right\rangle + \left\langle x^*,a_{i,k}^+\right\rangle\left\langle x^*,a_{i,k}^-\right\rangle \geqq 0 \iff$$

$$\left\langle x,a_{i,k}^+\right\rangle\left\langle x,a_{i,k}^-\right\rangle - \langle x,\alpha_{i,k}\rangle + \beta_{i,k} \geqq 0 \iff$$

$$x^T\left[a_{i,k}^+ {a_{i,k}^-}^T\right]x - \alpha_{i,k}^T x + \beta_{i,k} \geqq 0 \iff H_{i,k}(x) \geqq 0.$$

Therefore, the result follows from Theorem 3.4. ∎

From Corollary 3.5, if each quadratic function $H_{i,k}(x)$, $\{i,k\} \in I \times K$, is non-negative in the feasible set, then a feasible point satisfying (1), (2) and (3) is an Pareto optimal solution of (VQFP).

Let $\left[\bar{H}_{i,k}\right] \equiv \left[a_{i,k}^+ {a_{i,k}^-}^T\right] \in \mathbb{R}^{n\times n}$. Then $\beta_{i,k} = x^{*T}\left[\bar{H}_{i,k}\right]x^*$, and the non-negativity of the quadratic $H_{i,k}(x) = x^T\left[\bar{H}_{i,k}\right]x - \alpha_{i,k}^T x + \beta_{i,k}$ depends on each matrix $\left[\bar{H}_{i,k}\right]$ and each vector $\alpha_{i,k} \in \mathbb{R}^n$, where $\{i,k\} \in I \times K$ and $x^* \in S$. For example, the unconstrained (VQFP) requires that each matrix $\left[\bar{H}_{i,k}\right]$ be positive semidefinite and that $\beta_{i,k} \geqq \left(\alpha_{i,k} - w^T\left[\bar{H}_{i,k}\right]w\right)$, where $w$ is a solution of the system $2\left[\bar{H}_{i,k}\right]x = \alpha_{i,k}$.

**Corollary 3.6** *Let $x^*$ be a feasible point of (VQFP). Suppose that the constraint function $h_j$ is convex for each $j \in J$ and there exist vectors $\tau \in \mathbb{R}^m$, $\lambda \in \mathbb{R}^\ell$, such that (1), (2) and (3) are valid. If for each pair $\{i,k\} \in I \times K$, the matrix $\left[\bar{H}_{i,k}\right]$ is positive semidefinite and $\alpha_{i,k} = 0$ (see (12), (13) and (14)), then $x^* \in Eff(VQFP)$.*

*Proof* By hypothesis, for all $x \in S$ we have $H_{i,k}(x) \geqq 0$ for each pair $\{i,k\} \in I \times K$. Therefore, the result follows from Corollary 3.5. ∎

Given a pair $\{i,k\} \in I \times K$, writing each entry of the matrix $\left[\bar{H}_{i,k}\right] = \left(\bar{H}_{i,k}(r,s)\right)$ and each entry of the vector $\alpha_{i,k} = \left(\alpha_{i,k}(r)\right)$ according to the entries of the eigenvectors $p_i^k = \left(p_i^k(r)\right)$ and $q_i^k = \left(q_i^k(r)\right)$, where $r,s \in K$, we obtain for each pair $\{r,s\} \in K \times K$,

$$\begin{aligned} \bar{H}_{i,k}(r,s) &= a_{i,k}^+(r)a_{i,k}^-(s) = \\ &= \mu_k^{A_i}\, p_i^k(r)p_i^k(s) + \sqrt{\mu_k^{A_i}\mu_k^{B_i}\frac{f_i(x^*)}{g_i(x^*)}}\left(p_i^k(s)q_i^k(r) - p_i^k(r)q_i^k(s)\right) - \mu_k^{B_i}\frac{f_i(x^*)}{g_i(x^*)}\, q_i^k(r)q_i^k(s), \end{aligned} \tag{15}$$

$$\bar{H}_{i,k}(r,r) = \mu_k^{A_i}\left(p_i^k(r)\right)^2 - \mu_k^{B_i}\frac{f_i(x^*)}{g_i(x^*)}\left(q_i^k(r)\right)^2, \tag{16}$$

$$\alpha_{i,k}(r) = 2\mu_k^{A_i}\left\langle x^*, p_i^k\right\rangle p_i^k(r) - 2\mu_k^{B_i}\frac{f_i(x^*)}{g_i(x^*)}\left\langle x^*, q_i^k\right\rangle q_i^k(r). \tag{17}$$

We can draw some conclusions from (15), (16) and (17). For example, for a fixed pair $\{i,k\} \in I \times K$, the vector $\alpha_{i,k}$ is a linear combination of the eigenvectors $p_i^k$ and $q_i^k$. If $\mu_k^{A_i} = 0$, $\mu_k^{B_i} = 0$ or $f_i(x^*) = 0$, then $\left[\bar{H}_{i,k}\right]$ is a symmetric matrix. Moreover, if $\mu_k^{A_i}\mu_k^{B_i} f_i(x^*) < 0$, and there exists a pair $\{r,s\} \in K \times K$ such that $p_i^k(s)q_i^k(r) \neq p_i^k(r)q_i^k(s)$, then the matrix $\left[\bar{H}_{i,k}\right] \notin \mathbb{R}^{n\times n}$. In this case, if there exists $x \in S$ such that $H_{i,k}(x) \in \mathbb{C}\setminus\mathbb{R}$, $H_{i,k}(x) \geqq 0$ does not make sense. However, when (11) is required, it is possible to show that $H_{i,k}(x) \in \mathbb{C}\setminus\mathbb{R}$ is not possible.

The results of Theorem 3.1, 3.4 and its corollaries can be used in order to develop a method of searching for Pareto optimal solutions of (VQFP), and it might be useful to determine the termination criteria in the development of algorithms.

# 4 Duality

The matrix (10) defines a specific function, and by adding some assumptions about it, we obtain new results, such as, a relationship between the problem (VQFP) and a scalar problem associated with it, and the main duality theorems.

In the scalar optimization problem case, Dinkelbach [12] and Jagannathan [16] used a parametric approach that transforms the fractional optimization problem in a new scalar optimization problem. Similarly, we consider the following associated problem to (VQFP).

$$\begin{array}{lll}(\textbf{VQFP})_{x^*} & \textbf{Minimize} & f(x) - \dfrac{f(x^*)}{g(x^*)}g(x) = \left(f_1(x) - \dfrac{f_1(x^*)}{g_1(x^*)}g_1(x), \ldots, f_m(x) - \dfrac{f_m(x^*)}{g_m(x^*)}g_m(x)\right) \\ & \text{subject to} & h_j(x) \leqq 0 \quad j \in J, \\ & & x \in \Omega,\end{array}$$

where $\Omega \subseteq \mathbb{R}^n$, $f_i$, $g_i$, $i \in I$ and $h_j$, $j \in J$ are defined in (VQFP), and $x^* \in S$.

Using assumptions of generalized convexity, Osuna-Gómez et al. [30] presented the problem $(\text{VFP})_{x^*}$ and obtained necessary and sufficient conditions for weakly Pareto optimality and main duality theorems. The results presented in [12, 16, 30] considered each objective function as $f_i(x) - \alpha_i g_i(x)$, $i \in I$, and they studied the properties of the parameter $\alpha_i \in \mathbb{R}$. Following the ideas presented by Osuna-Gómez et al. [30], we obtain new results by considering directly $\alpha_i \equiv \frac{f_i(x^*)}{g_i(x^*)}$, $i \in I$, where $x^* \in Eff(VQFP)$. However, by imposing hypothesis on the linear combination of matrices $\left[A_i - \frac{f_i(x^*)}{g_i(x^*)}B_i\right]$, $i \in I$ and $x^* \in S$, we consider Pareto optimal solutions rather than weakly Pareto optimal solutions.

To characterize the solutions of the problems (VOP), Geoffrion [13] used the solutions of the associated scalar problems. Similarly, we consider the following *weighted scalar problem* associated to the problem $(\text{VQFP})_{x^*}$.

$$\begin{array}{lll}(\textbf{VQFP})_{x^*}^{w} & \textbf{Minimize} & \sum\limits_{i=1}^{m} w_i\left(f_i(x) - \frac{f_i(x^*)}{g_i(x^*)}g_i(x)\right) \\ & \text{subject to} & h_j(x) \leqq 0 \quad j \in J, \\ & & x \in \Omega,\end{array}$$

where $\Omega \subseteq \mathbb{R}^n$, $f_i$, $g_i$, $i \in I$ and $h_j$, $j \in J$ are defined in (VQFP), $x^* \in S$ and $w = (w_1, \ldots, w_m)^T \in \mathbb{R}^m$, $w > 0$.

## 4.1 The relationship between the associated problems

The next theorem and its proof are similar to Lemma 1.1 from [30], when Pareto optimal solutions (not necessarily weak) are considered.

**Theorem 4.1** *$x^* \in Eff(VQFP)$ if and only if $x^* \in Eff(VQFP)_{x^*}$.*

*Proof* See Lemma 1.1 in [30], considering "$\leq$" instead of "$<$". ∎

In Section 3, we define the matrix $\widehat{F}(x^*) = \sum_{i=1}^{m} \left[ \frac{\tau_i}{g_i(x^*)} A_i - \frac{\tau_i f_i(x^*)}{[g_i(x^*)]^2} B_i \right]$, where $x^* \in S$ and $\tau \in \mathbb{R}^m$, $\tau > 0$. Let us define now the set $\mathcal{W} = \{w \in \mathbb{R}^m \mid w > 0\}$, the function $F : \mathcal{W} \times S \to \mathbb{R}^{n \times n}$ given by

$$F(w,x) \equiv \sum_{i=1}^{m} w_i \left[ A_i - \frac{f_i(x)}{g_i(x)} B_i \right],$$

and for each $i \in I$, the functions $F_i : S \to \mathbb{R}^{n \times n}$ given by

$$F_i(x) \equiv A_i - \frac{f_i(x)}{g_i(x)} B_i.$$

Then, we have $\widehat{F}(x^*) = F\left(\frac{\tau}{g(x^*)}, x^*\right)$, where $\frac{\tau}{g(x^*)} = \left(\frac{\tau_1}{g_1(x^*)}, \ldots, \frac{\tau_m}{g_m(x^*)}\right)^T \in \mathcal{W}$, $F(w,x) = \sum_{i=1}^{m} w_i F_i(x)$ and we can establish some relations among the associated problems (VQFP), $(\text{VQFP})_{x^*}$ and $(\text{VQFP})_{x^*}^{w}$.

**Theorem 4.2** *If $x^*$ is a optimal solution of the weighted scalar problem $(VQFP)_{x^*}^{w}$, then $x^* \in Eff(VQFP)$.*

*Proof* Suppose that $x^* \notin Eff(VQFP)$, then there exists another point $x \in S$ such that

$$\frac{f_i(x)}{g_i(x)} \leq \frac{f_i(x^*)}{g_i(x^*)} \implies f_i(x) - \frac{f_i(x^*)}{g_i(x^*)} g_i(x) \leq 0 \implies \sum_{i=1}^{m} w_i \left( f_i(x) - \frac{f_i(x^*)}{g_i(x^*)} g_i(x) \right) < 0$$

$$\implies \sum_{i=1}^{m} w_i \left( f_i(x) - \frac{f_i(x^*)}{g_i(x^*)} g_i(x) \right) < \sum_{i=1}^{m} w_i \left( f_i(x^*) - \frac{f_i(x^*)}{g_i(x^*)} g_i(x^*) \right).$$

This contradicts the minimality of $x^*$ in $(\text{VQFP})_{x^*}^{w}$. ∎

**Lemma 4.1** *Let $x^* \in Eff(VQFP)$. Suppose that the constraint qualification (GGCQ) is satisfied at $x^*$, then there exist vectors $\tau^* > 0$ and $\lambda^* \geqq 0$ such that*

$$\sum_{i=1}^{m} \tau_i^* \left( \nabla f_i(x^*) - \frac{f_i(x^*)}{g_i(x^*)} \nabla g_i(x^*) \right) + \sum_{j=1}^{\ell} \lambda_j^* \nabla h_j(x^*) = 0,$$

$$\sum_{j=1}^{\ell} \lambda_j^* h_j(x^*) = 0.$$

*Proof* Let $x^* \in S$, $\mu \in \mathbb{R}^m$, $\mu > 0$ and $\tau_i = \frac{\mu_i}{g_i(x^*)} > 0$, $i \in I$. Then

$$\begin{aligned}
\sum_{i=1}^{m} \mu_i \nabla \frac{f_i(x^*)}{g_i(x^*)} &= \sum_{i=1}^{m} \tau_i g_i(x^*) \, \nabla \frac{f_i(x^*)}{g_i(x^*)} \\
&= \sum_{i=1}^{m} \tau_i g_i(x^*) \left( \frac{\nabla f_i(x^*) g_i(x^*) - \nabla g_i(x^*) f_i(x^*)}{[g_i(x^*)]^2} \right) \\
&= \sum_{i=1}^{m} \tau_i \left( \nabla f_i(x^*) - \frac{f_i(x^*)}{g_i(x^*)} \nabla g_i(x^*) \right), \qquad (18)
\end{aligned}$$

and if $x^* \in \mathit{Eff}(VQFP)$, by Lemma 2.1, there exist $\mu^* > 0$ and $\lambda^* \geqq 0$ such that $(x^*, \mu^*, \lambda^*)$ is a critical point, in the $KKT$ sense, of the problem (VQFP). That is,

$$\begin{aligned}
&\sum_{i=1}^{m} \mu_i^* \nabla \frac{f_i(x^*)}{g_i(x^*)} + \sum_{j=1}^{\ell} \lambda_j^* \nabla h_j(x^*) = 0, \\
&\sum_{j=1}^{\ell} \lambda_j^* h_j(x^*) = 0.
\end{aligned}$$

From (18), there exist $\tau^* > 0$, $\tau_i^* = \frac{\mu_i^*}{g_i(x^*)} > 0$, $i \in I$, and $\lambda^* \geqq 0$ such that

$$\begin{aligned}
&\sum_{i=1}^{m} \tau_i^* \left( \nabla f_i(x^*) - \frac{f_i(x^*)}{g_i(x^*)} \nabla g_i(x^*) \right) + \sum_{j=1}^{\ell} \lambda_j^* \nabla h_j(x^*) = 0, \\
&\sum_{j=1}^{\ell} \lambda_j^* h_j(x^*) = 0.
\end{aligned}$$

Therefore, the result is valid. ■

**Lemma 4.2** *Let $x^* \in S$. If there exists $w \in \mathcal{W}$, such that the matrix $F(w, x^*)$ is positive semidefinite, then the objective function of $(VQFP)_{x^*}^{w}$ is convex.*

*Proof* Given $x_1, x_2 \in S$, we have for each $i \in I$,

$$\begin{aligned}
f_i(x_1) - f_i(x_2) &= (x_1 - x_2)^T A_i (x_1 - x_2) + \nabla f_i(x_2)^T (x_1 - x_2), \\
g_i(x_1) - g_i(x_2) &= (x_1 - x_2)^T B_i (x_1 - x_2) + \nabla g_i(x_2)^T (x_1 - x_2).
\end{aligned}$$

Hence, for each objective function of $(VQFP)_{x^*}$, we have

$$\left( f_i(x_1) - \frac{f_i(x^*)}{g_i(x^*)} g_i(x_1) \right) - \left( f_i(x_2) - \frac{f_i(x^*)}{g_i(x^*)} g_i(x_2) \right) =$$

$$= \left( f_i(x_1) - f_i(x_2) \right) - \frac{f_i(x^*)}{g_i(x^*)} \left( g_i(x_1) - g_i(x_2) \right)$$

$$= \left\{(x_1-x_2)^T A_i (x_1-x_2) + \nabla f_i(x_2)^T (x_1-x_2)\right\} +$$

$$-\frac{f_i(x^*)}{g_i(x^*)}\left\{(x_1-x_2)^T B_i (x_1-x_2) + \nabla g_i(x_2)^T (x_1-x_2)\right\}$$

$$= (x_1-x_2)^T \left[A_i - \frac{f_i(x^*)}{g_i(x^*)} B_i\right](x_1-x_2) + \left(\nabla f_i(x_2) - \frac{f_i(x^*)}{g_i(x^*)}\nabla g_i(x_2)\right)^T (x_1-x_2). \quad (19)$$

If there exists $w \in \mathcal{W}$ such that the matrix $F(w,x^*)$ is positive semidefinite, then

$$\sum_{i=1}^m w_i \left(f_i(x_1) - \frac{f_i(x^*)}{g_i(x^*)} g_i(x_1)\right) - \sum_{i=1}^m w_i \left(f_i(x_2) - \frac{f_i(x^*)}{g_i(x^*)} g_i(x_2)\right) =$$

$$= (x_1-x_2)^T \sum_{i=1}^m w_i \left[A_i - \frac{f_i(x^*)}{g_i(x^*)} B_i\right](x_1-x_2) + \sum_{i=1}^m w_i \left(\nabla f_i(x_2) - \frac{f_i(x^*)}{g_i(x^*)}\nabla g_i(x_2)\right)^T (x_1-x_2)$$

$$\geq \left(\sum_{i=1}^m w_i \left(\nabla f_i(x_2) - \frac{f_i(x^*)}{g_i(x^*)}\nabla g_i(x_2)\right)\right)^T (x_1-x_2).$$

Therefore, the objective function of $(VQFP)_{x^*}^w$ is convex. ■

Note that the hypothesis of semidefiniteness on the matrix $F(w,x^*)$ or on the matrices $F_i(x^*)$, $i \in I$, $x^* \in S$, is punctual. However, the next example, we draw a situation in which for all $x \in S$ and $i \in I$, we have $y^T F_i(x) y \geqq 0$, for all $y \in S$, and then $y^T F(w,x) y \geqq 0$, for all $y \in S$.

*Example* Consider the problem (VQFP), where $S = [-2,2]$ and for all $x \in S$

$$\frac{f_1(x)}{g_1(x)} = \frac{x-2}{x^2+2}, \quad \frac{f_2(x)}{g_2(x)} = \frac{2x^2-x-1}{x^2+1}, \quad \frac{f_3(x)}{g_3(x)} = \frac{-2x^2-2x-5}{x^2+x+1}.$$

For all these functions, we obtain for all $y \in S$

$$\begin{aligned}
y^T \left[A_1 - \frac{f_1(x)}{g_1(x)} B_1\right] y &= y^2 \left(0 - \frac{x-2}{x^2+2}\,1\right) = \frac{y^2(2-x)}{x^2+2} \geqq 0,\\
y^T \left[A_2 - \frac{f_2(x)}{g_2(x)} B_2\right] y &= y^2 \left(2 - \frac{2x^2-x-1}{x^2+1}\,1\right) = \frac{y^2(x+3)}{x^2+1} \geqq 0,\\
y^T \left[A_3 - \frac{f_3(x)}{g_3(x)} B_3\right] y &= y^2 \left(-2 - \frac{-2x^2-2x-5}{x^2+x+1}\,1\right) = \frac{3y^2}{x^2+x+1} \geqq 0.
\end{aligned}$$

Therefore, for this example, each point $x^*$ satisfying (1), (2) and (3) is Pareto optimal. For example, for $(\tau_1,\tau_2,\tau_3)^T = (0.5,1,0.25)^T$ and $(\lambda_1,\lambda_2)^T = (1,1)$, we have that $x^* = 0$ is Pareto optimal solution. Likewise, for $(\tau_1,\tau_2,\tau_3)^T \approx (0.62,1,0.89)^T$ and $(\lambda_1,\lambda_2)^T = (0,0)^T$, we have that $x^* = -0.25$ is Pareto optimal solution. ■

Theorem 4.1 shows an equivalence between the associated problems (VQFP) and (VQFP)$_{x^*}$. In the next theorem shows a relation between the problems (VQFP)$_{x^*}$ and (VQFP)$^w_{x^*}$, then it provides a converse to the Theorem 4.2.

**Theorem 4.3** *Let $x^* \in Eff(VQFP)_{x^*}$. Suppose that the constraint qualification (GGCQ) is satisfied at $x^*$, and the constraint function $h_j$ is convex for each $j \in J$. Then there exists $w \in \mathcal{W}$ such that if the matrix $F(w,x^*)$ is positive semidefinite, then $x^*$ is the optimal solution for the weighted scalar problem (VQFP)$^w_{x^*}$.*

*Proof* If $x^* \in Eff(VQFP)_{x^*}$ and satisfies (GGCQ), by Lemma 4.1, there exist $w > 0$ and $\lambda \geqq 0$, such that

$$\begin{aligned}
&\sum_{i=1}^{m} w_i^* \left( \nabla f_i(x^*) - \frac{f_i(x^*)}{g_i(x^*)} \nabla g_i(x^*) \right) + \sum_{j=1}^{\ell} \lambda_j^* \nabla h_j(x^*) = 0, \\
&\sum_{j=1}^{\ell} \lambda_j^* h_j(x^*) = 0.
\end{aligned}$$

Therefore, $x^*$ is a critical point of the weighted scalar problem (VQFP)$^w_{x^*}$, and since $F(w,x^*)$ is positive semidefinite, by Lemma 4.2, the objective function of (VQFP)$^w_{x^*}$ is convex. Since for each $j \in J$ the constraint function $h_j$ is convex, it follows that $x^*$ is an optimal solution for (VQFP)$^w_{x^*}$. ■

## 4.2 Duality theorems

For a given mathematical optimization problem there are many types of duality. Two well-known duals are the Wolfe dual [38] and the Mond-Weir dual [29]. In this work, we consider the primal problem (VQFP) and discuss the Mond-Weir dual problem, but we use the associated problem (VQFP)$_{x^*}$ to generate the constraint set of the dual problem. Let us consider the following vector quadratic fractional dual optimization problem (VQFD).

$$\begin{array}{lll}
\textbf{(VQFD)} & \textbf{Maximize} & \frac{f(u)}{g(u)} = \left( \frac{f_1(u)}{g_1(u)}, \ldots, \frac{f_m(u)}{g_m(u)} \right) \\
 & \text{subject to} & \sum_{i=1}^{m} \tau_i \left( \nabla f_i(u) - \frac{f(u)}{g(u)} \nabla g_i(u) \right) + \sum_{j=1}^{\ell} \lambda_j \nabla h_j(u) = 0, \\
 & & \sum_{j=1}^{\ell} \lambda_j h_j(u) \geqq 0, \\
 & & \tau > 0, \quad \lambda \geqq 0, \quad \sum_{j=1}^{\ell} \lambda_j = 1, \\
 & & u \in S,
\end{array}$$

where $f_i$ and $g_i$, $i \in I$ are the same quadratic functions defined on (VQFP), and its feasible set we denote by $Y$.

**Theorem 4.4** *(Weak duality) Let $x \in S$ and $(u,\tau,\lambda) \in Y$. If $F(\tau,u)$ is positive semidefinite, and the constraint function $h_j$ is convex for each $j \in J$, then*

$$\frac{f(x)}{g(x)} \nleq \frac{f(u)}{g(u)}.$$

*Proof* If there are $x \in S$ and $(u,\tau,\lambda) \in Y$ such that $\frac{f(x)}{g(x)} \leq \frac{f(u)}{g(u)}$, then

$$\frac{f(x)}{g(x)} \leq \frac{f(u)}{g(u)} \implies f(x) - \frac{f(u)}{g(u)}g(x) \leq 0 \implies \sum_{i=1}^{m} \tau_i \left( f_i(x) - \frac{f_i(u)}{g_i(u)}g_i(x) \right) < 0.$$

Since $x \in S$, then $\sum_{j=1}^{\ell} \lambda_j h_j(x) \leqq 0$, and $\sum_{j=1}^{\ell} \lambda_j h_j(u) \geqq 0$ implies that

$$\sum_{i=1}^{m} \tau_i \left( f_i(x) - \frac{f_i(u)}{g_i(u)}g_i(x) \right) + \sum_{j=1}^{\ell} \lambda_j h_j(x) < \sum_{i=1}^{m} \tau_i \left( f_i(u) - \frac{f_i(u)}{g_i(u)}g_i(u) \right) + \sum_{j=1}^{\ell} \lambda_j h_j(u).$$

Once $F(\tau,u)$ is positive semidefinite and each constraint function $h_j$ is convex, we can use Lemma 4.2 to conclude that the objective function of (VQFP)$^{\tau}_{x^*}$ is convex, and

$$\begin{aligned}
0 &> \sum_{i=1}^{m} \tau_i \left[ \left( f_i(x) - \frac{f_i(u)}{g_i(u)}g_i(x) \right) - \left( f_i(u) - \frac{f_i(u)}{g_i(u)}g_i(u) \right) \right] + \sum_{j=1}^{\ell} \lambda_j \left( h_j(x) - h_j(u) \right) \\
&\geqq \sum_{i=1}^{m} \tau_i \left( \nabla f_i(u) - \frac{f_i(u)}{g_i(u)} \nabla g_i(u) \right)^T (x-u) + \sum_{j=1}^{\ell} \lambda_j \nabla h_j(u)^T (x-u) \\
&= (x-u)^T \left( \sum_{i=1}^{m} \tau_i \left( \nabla f_i(u) - \frac{f_i(u)}{g_i(u)} \nabla g_i(u) \right) + \sum_{j=1}^{\ell} \lambda_j \nabla h_j(u) \right) = 0,
\end{aligned}$$

which is a contradiction. ∎

**Theorem 4.5** *(Strong duality) Let $x^* \in Eff(VQFP)$. Suppose that (GGCQ) holds at $x^*$, then there exists $(\tau^*,\lambda^*)$ such that $(x^*,\tau^*,\lambda^*)$ is feasible for (VQFD) and the values of the objective function of (VQFP) and (VQFD) are equal. Moreover, if $F(\tau^*,x^*)$ is positive semidefinite, and the constraint function $h_j$ is convex for each $j \in J$, then $(x^*,\tau^*,\lambda^*) \in Eff(VQFD)$.*

*Proof* If $x^* \in Eff(VQFP)$, by Lemma 4.1, there are $\tau^* > 0$ and $\lambda^* \geqq 0$ such that $(x^*,\tau^*,\lambda^*)$ satisfies

$$\begin{aligned}
&\sum_{i=1}^{m} \tau_i^* \left( \nabla f_i(x^*) - \frac{f_i(x^*)}{g_i(x^*)} \nabla g_i(x^*) \right) + \sum_{j=1}^{\ell} \lambda_j^* \nabla h_j(x^*) = 0, \\
&\sum_{j=1}^{\ell} \lambda_j^* h_j(x^*) = 0.
\end{aligned}$$

Then $(x^*,\tau^*,\lambda^*) \in Y$ and the values of the objective functions of (VQFP) and (VQFD) are equal. Moreover, if $F(\tau^*,x^*)$ is positive semidefinite, each constraint function $h_j$ is convex, and $(x^*,\tau^*,\lambda^*) \notin Eff(VQFD)$, then there exists another point $(u,\tau,\lambda) \in Y$ such that

$$\frac{f(u)}{g(u)} \geq \frac{f(x^*)}{g(x^*)},$$

contradicting the weak duality. ∎

**Theorem 4.6** *(Converse duality) Let $(u^*,\tau^*,\lambda^*) \in Y$ and $u^*$ be feasible point of the primal problem (VQFP). If $F(\tau^*,u^*)$ is positive semidefinite, and the constraint function $h_j$ is convex for each $j \in J$, then $u^* \in Eff(VQFP)$.*

*Proof* If $(u^*,\tau^*,\lambda^*) \in Y$ and $u^* \in S$, then $\sum_{j=1}^{\ell} \lambda_j^* h_j(u^*) \geqq 0$, $\sum_{j=1}^{\ell} \lambda_j^* h_j(u^*) \leqq 0$ and

$$\sum_{i=1}^{m} \tau_i^* \left( \nabla f_i(u^*) - \frac{f_i(u^*)}{g_i(u^*)} \nabla g_i(u^*) \right) + \sum_{j=1}^{\ell} \lambda_j^* \nabla h_j(u^*) = 0,$$
$$\sum_{j=1}^{\ell} \lambda_j^* h_j(u^*) = 0.$$

Therefore, $u^*$ is a critical point for the weighted scalar problem $(\text{VQFP})_{u^*}^{\tau^*}$. Since $F(\tau^*,u^*)$ is positive semidefinite, by Lemma 4.2, the objective function of $(\text{VQFP})_{u^*}^{\tau^*}$ is convex. Moreover, if each constraint function $h_j$ is convex, $j \in J$, then $u^*$ is an optimal solution of $(\text{VQFP})_{u^*}^{\tau^*}$. Thus, by Theorem 4.2, we have $u^* \in Eff(VQFP)$. ■

We can obtain a second type of converse duality theorem requiring more of the matrix function $F$. Specifically, there must be vectors $(w,x) \in \mathcal{W} \times S$ such that $F(w,x)$ is positive definite, that is, $y^T F(w,x) y > 0$, $\forall y \in \mathbb{R}^n$ and $y \neq 0$.

**Theorem 4.7** *(Strict converse duality) Let $x^* \in S$ and $(u^*,\tau^*,\lambda^*) \in Y$ such that*

$$\sum_{i=1}^{m} \tau_i^* \frac{f_i(x^*)}{g_i(x^*)} = \sum_{i=1}^{m} \tau_i^* \frac{f_i(u^*)}{g_i(u^*)}. \tag{20}$$

*If the matrix $F\left(\frac{\tau^*}{g(x^*)},u^*\right)$ is positive definite and the constraint function $h_j$ is convex for each $j \in J$, then $x^* = u^*$.*

*Proof* Suppose $x^* \neq u^*$. Since $x^* \in S$ and $(u^*,\tau^*,\lambda^*) \in Y$, then $-\sum_{j=1}^{\ell} \lambda_j^* h_j(u^*) \leqq 0$ and $\sum_{j=1}^{\ell} \lambda_j^* h_j(x^*) \leqq 0$. If each constraint function $h_j$ is convex, $j \in J$, we obtain

$$\begin{aligned}
0 \quad &\geqq \quad \sum_{j=1}^{\ell} \lambda_j^* h_j(x^*) - \sum_{j=1}^{\ell} \lambda_j^* h_j(u^*) \geqq \left( \sum_{j=1}^{\ell} \lambda_j^* \nabla h_j(u^*) \right)^T (x^* - u^*) = \\
&= \quad \left( -\sum_{i=1}^{m} \tau_i^* \left( \nabla f_i(u^*) - \frac{f_i(u^*)}{g_i(u^*)} \nabla g_i(u^*) \right) \right)^T (x^* - u^*) \quad \Longrightarrow \\
&\Longrightarrow \quad \sum_{i=1}^{m} \tau_i^* \left( \nabla f_i(u^*) - \frac{f_i(u^*)}{g_i(u^*)} \nabla g_i(u^*) \right)^T (x^* - u^*) \geqq 0.
\end{aligned}$$

Using the proof of Theorem 3.1, given $u^* \in S$ and $i \in I$, for all $x \in S$, we have

$$\frac{f_i(x)}{g_i(x)} - \frac{f_i(u^*)}{g_i(u^*)} = \frac{1}{g_i(x)} \left\{ (x - u^*)^T \left[ A_i - \frac{f_i(u^*)}{g_i(u^*)} B_i \right] (x - u^*) \right\} +$$

$$+\frac{1}{g_i(x)}\left\{\left(\nabla f_i(u^*)-\frac{f_i(u^*)}{g_i(u^*)}\nabla g_i(u^*)\right)^T(x-u^*)\right\}.$$

Therefore, for $(\tau^*,x^*)\in\mathcal{W}\times S$ we obtain

$$\sum_{i=1}^m\tau_i^*\left(\frac{f_i(x^*)}{g_i(x^*)}-\frac{f_i(u^*)}{g_i(u^*)}\right)=(x^*-u^*)^T\left[F\left(\frac{\tau^*}{g(x^*)},u^*\right)\right](x^*-u^*)+$$
$$+\sum_{i=1}^m\frac{\tau_i^*}{g_i(x^*)}\left(\nabla f_i(u^*)-\frac{f_i(u^*)}{g_i(u^*)}\nabla g_i(u^*)\right)^T(x^*-u^*),$$

and since $F\left(\frac{\tau^*}{g(x^*)},u^*\right)$ is positive definite and $x^*\neq u^*$, then by (20)

$$\begin{aligned}0 &= \sum_{i=1}^m\tau_i^*\left(\frac{f_i(x^*)}{g_i(x^*)}-\frac{f_i(u^*)}{g_i(u^*)}\right)>\sum_{i=1}^m\frac{\tau_i^*}{g_i(x^*)}\left(\nabla f_i(u^*)-\frac{f_i(u^*)}{g_i(u^*)}\nabla g_i(u^*)\right)^T(x^*-u^*)\\ &\geqq \min_{i\in I}\left\{\frac{1}{g_i(x^*)}\right\}\sum_{i=1}^m\tau_i^*\left(\nabla f_i(u^*)-\frac{f_i(u^*)}{g_i(u^*)}\nabla g_i(u^*)\right)^T(x^*-u^*)\geqq 0,\end{aligned}$$

which is a contradiction. ■

# 5 Conclusions

The main contribution of this work is the development of Pareto optimality conditions for a particular vector optimization problem, where each objective function consists of a ratio of two quadratic functions with convexity being only assumed on the constraint set. We took advantage of the diagonalization of Hessian matrices. We have shown the relationship between the particular problem and two problems associated with it, and we use some assumptions of the linear combination of Hessian matrices to show the main duality theorems. For the particular problem, the results presented in this work might be useful to determine the termination criteria in the development of algorithms, and new extensions can be established to more general vector optimization problems, in which algorithms based on quadratic approximations are used locally. In future work we plan to develop algorithms using the concepts presented here.

# Acknowledgements

We are indebted to the anonymous reviewers for their helpful comments. The first author was supported by Coordination for the Improvement of Higher Level Personnel of Brazil (CAPES). The second author was partially supported by Spains Ministry of Science and Technology under grant MTM2007-63432. The third and fourth authors were partially supported by National Council for Scientific and Technological Development of Brazil (CNPq) and Foundation for Research Support of the State of São Paulo (FAPESP).